\documentstyle[12pt,leqno]{article}
\newcommand{\se}[1]{{\section{#1}} {\setcounter{equation}{0}}}
\newtheorem{th}{Theorem}[section]
\newtheorem{lm}{Lemma}[section]
\newtheorem{prop}{Proposition}[section]
\newtheorem{de}{Definition}[section]
\newtheorem{co}{Corollary}[section]
\newtheorem{re}{Remark}[section]
\def\k{{K\"{a}hler }}
\def\ke{{K\"{a}hler-Einstein }}

\begin{document}
\hbadness=10000
\title{Canonical coordinates and Bergmann metrics}
\author{Wei-Dong Ruan\\
Department of mathematics\\
Columbia University\\
New York, NY 10027\\}
\footnotetext{Partially supported by Sonderforschungsbereich 237.}
\maketitle
\begin{abstract}
In this paper we will discuss local coordinates canonically corresponding to a \k metric. We will also discuss the $C^\infty$ convergence of Bergmann metrics following Tian's result on $C^2$ convergence of Bergmann metrics. At the end we present an interesting characterization of ample line bundle that could be useful in arithmetic geometry.
\end{abstract}
\newpage
\tableofcontents
\newpage
\section{Introduction and Background}
For real analytic \k manifold, there is a canonically defined local coordinate in a neighborhood of any point, which depends antiholomorphically on the point. Since the coordinates canonically correspond to \k metrics, we will call them K-coordinates, representing \k coordinates. These coordinates first come to our attention from the literature in mirror symmetry. It is given a nice geometric interpretation by Kontsevich. We will discuss this interpretation, write it in terms of local coordinates, and see its geometric implications.\\\\
K-coordinates actually were quite well known in mathematics literature. To author's knowledge, they first appeared in [B]. Given the ``canonical'' nature of these coordinates, it is quite surprising to me that they were not widely used in the literature, although author first saw a finite order approximation of it being used in Tian's proof of $C^2$ convergence of Bergmann metrics. It is very interesting, at least to the author, to see how useful are these coordinates in handling problems in \k geometry and other fields. We will illustrate their use through one example: the proof of $C^\infty$ convergence of Bergmann metrics.\\\\
It was first conjectured by Yau that Bergmann metrics induced by a \ke metric should converge to original \ke metric. Later, Tian proved that Bergmann metrics induced by any \k metric of integral \k class will converge to that \k metric in $C^2$ norm. The problem is roughly as follows. For $(M, g)$ a \k manifold with the \k class $\omega_g \in H^2(M, {\bf Z})$, one can realize $\omega_g$ as first Chern class of an ample line bundle $L$. Multiples $L^m$ will induce certain canonical embeddings to Projective spaces by all its sections for $m$ large.
\[
\varphi_m: M \longrightarrow {\bf CP^{N_m}}
\]
pull backs of Fubini-Study metrics are usually called Bergmann metrics, they look like a nice collection of metrics canonically associated to the \k manifold, very much like polynomials. It then sound natural to expect that a suitably chosen sequence of Bergmann metrics will converge to the original \k metric. The precise definitions and statements will be given in section 3.\\\\
The starting point for Tian's proof is the construction of a sequence of ``peak sections'', which, when used correctly, will reduce the problem to certain kind of computation. He then handle these computations by brute force, for which he already run into quite complicated computation in proving the $C^2$ convergence. Higher order convergences do not seem manageable with the same way. Although he believe that the $C^\infty$ convergence should be true.\\\\
In our way of understanding the K-coordinates, we find his proof very natural in K-coordinate point of view, which enable us to choose a more canonical factor in a lemma to make it simpler. This turns out to be a crucial modification, which makes the resulting computation problem more tackable, at least in theory, it turns into a pure combinatoric problem, which, in principle, should be solvable. But it turns out the corresponding combinatoric problem is very complicated to solve, if not impossible. (I have more respect and appreciation to combinatoric problems now than before.) Although we were able to solve it to a great extent, even deduced an algorithm to prove convergence to arbitrary order, modulo computations that get more and more complicated. We are still not able to solve it by purely combinatoric computation.\\\\
Finally, it turns out, there is a very simple and elegant geometric proof, which even greatly simplifies the Tian's proof for the $C^2$ convergence. It reduce the needed computation to very minimal. But it do not seems to uncover the combinatoric mystery of the problem, which involve quite a lot of magical vanishing, which is quite interesting in its own.\\\\
Our paper is organized in the following way. In section 2, we will discuss the construction of K-coordinates, especially Kontsevich's interpretation. In section 3, we will introduce the convergence problem of Bergmann metrics, Tian's result, and carry out some delicate argument to reduce the problem to a combinatoric one. In section 4, We will demonstrate the combinatorics involved, and set up notation and combinatoric framework. In section 5, we try to sort out the complicated combinatorics to a great extent. While these problems is very interesting and challenging, they are not absolutely necessary, and not very easy to get through. Therefore can be ignored without serious lose. In section 6, we will present the simple geometric proof. In the last section, we will present some interesting application of our results.\\\\
\se{Construction of K-coordinates}
Before getting into the construction of K-coordinates, let's first review some facts from symplectic geometry. Given a symplectic manifold $(M, \omega)$ with a Lagrangian foliation $\pi: M \longrightarrow Q$, where each fiber of $\pi$ is a Lagrangian submanifold. On the tangent bundle of each Lagrangian submanifold, one can canonically define a flat connection by:
\[
\nabla_u v =i_\omega^{-1} {\cal L}_u i_\omega v
\]
where $u, v$ are vector fields along Lagrangian foliation, and $i_\omega: T_M\rightarrow T_M^*$ is defined as $i_\omega v = i(v)\omega$.\\\\
Assume $w$ is a vector field along Lagrangian foliation, we need to check $<{\cal L}_u i_\omega v, w>=0$, and ``$\nabla$'' as defined above is a connection.
\[
<{\cal L}_u i_\omega v, w>=u(\omega(v, w))-\omega (v, [u, w])=0
\]
and
\[
{\cal L}_u i_\omega (a v)=u(a)i_\omega v + a{\cal L}_u i_\omega v
\]
\[
{\cal L}_{(au)} i_\omega v = i(au) d(i_\omega v) + d (i(au)(i_\omega v))=a{\cal L}_u i_\omega v \\
\]
Recall that:
\[
{\cal L}_u {\cal L}_v \alpha - {\cal L}_v {\cal L}_u \alpha ={\cal L}_ {[u, v]} \alpha
\]
we get the flatness of the connection, therefore define a canonical affine structure on each Lagrangian fiber.\\\\
We can also see this in terms of local coordinates, let $q$ be the coordinates of $Q$, and $(p,q)$ be coordinates of $M$. In these coordinate
\[
\omega = b_{ij} dp^i\wedge dq^j + a_{ij} dq^i\wedge dq^j.
\]
The closeness of $\omega$ implies that $\partial(b_{ij})/\partial p^k = \partial(b_{kj})/\partial p^i$, which imply that
\[
\omega =  df^i\wedge dq^i + a_{ij} dq^i\wedge dq^j
\]
for some $f^i$. (Here $a_{ij}$ actually changed due to $q$ derivatives of $f^i$.) Now change the coordinates as $p^i=f^i$, we get
\[
\omega =  dp^i\wedge dq^i + a_{ij} dq^i\wedge dq^j
\]
here $a_{ij}$ depend only on $q$. $p$ so defined, gives us the affine structure on every Lagrangian leaf as induced by the connection. The ambiguities on determining $p$ are exactly affine transformations depend on $q$.\\\\
We can further simplify the expression of $\omega$ by noticing that $a_{ij} dq^i\wedge dq^j$ is closed, so\\
\[
a_{ij} dq^i\wedge dq^j = d(a_i dq^i) = d a_i\wedge dq^i
\]
and a transformation that depends on $q$: $p^i \longrightarrow p^i + a_i(q)$, will reduce $\omega$ to the standard form\\
\[
\omega =  dp^i\wedge dq^i
\]\\
Without lose of generality, we may fix the coordinate $q$, then the ambiguities of the coordinate $p$ are translations of following type\\
\[
p^i \longrightarrow p^i + \frac{\partial f(q)}{\partial q^i}
\]
Here $f(q)$ is a function that only depends on $q$.\\\\ 
Now, we recall the construction by Kontsevich. Given a \k manifold $(M, \omega)$, where $\omega$ is real analytic. Denote by $\overline{M}$ the same manifold $M$ endowed with the complex structure conjugate to the original one. The diagonal submanifold $M^{diag}$ of $M\times\overline{M}$ is totally real. Hence the differential form $\omega$ on $M^{diag}$ has the analytic continuation to the holomorphic form $\omega^{\bf C}$ in a neighborhood $U$ of $M^{diag}$. Thus $U$ is a complex symplectic manifold. It is easy to see that the submanifolds $M \times \{\overline{m}\}\cap U$, where $\overline{m}\in \overline{M}$, are Lagrangian. It means that we have a Lagrangian foliation of U. By our above discussion, leaves of such a foliation carry a natural flat affine structure. Hence this gives the construction of holomorphic affine structures on open subsets of $M$ depending antiholomorphically on points of $M$.\\\\ 
Let $(z, w)$ be local coordinates of $M\times\overline{M}$, then $i: M \rightarrow M^{diag} \hookrightarrow M\times\overline{M}$ can be written as: $i(z) = (z, \overline{z})$. Assume
\[
\omega = g_{ij}(z, \overline{z})dz^i \wedge d{\overline{z}}^j
\]
then
\[
\omega^{\bf C}= g_{ij}(z, w)dz^i \wedge dw^j
\]
For any $m \in M$, let $U_m = i^{-1}(M \times \{\overline{m}\}\cap U)$, and $\phi_m : U_m \longrightarrow {\bf C^n}$ is a flat coordinate with respect to m, as defined above, such that $\phi_m(m)=0$ and tangent map at $m$ is an isometry with respect to standard metric on ${\bf C^n}$.\\
\begin{re}
Maps $\phi_m$ are actually parameterized by elements on the total space of the canonical principle $U(n)$ bundle on $M$ with respect to the \k metric.\\
\end{re}
If we take the coordinates $z=\phi_m$, since $z$ is already flat coordinate, we have $g_{ij}(z, 0) = \delta_{ij}$, i.e. in the Taylor expansion of $g_{ij}$, all the $(0, l), \ (l, 0)$ terms vanish. This is equivalent to saying that in each $ U_m$ we can get a canonical \k potential $K_m$ , which depend on $m$, such that all the $(0, l), \ (1, l), \ (l, 1), \ (l, 0)$ terms in the Taylor expansion of $K_m(z)$ vanish except for $(1,1)$ term which equals to $|z|^2$. Locally 
\[
K_m(z) = |z|^2 + R_{i\bar{j}k\bar{l}}z^i{\bar{z}}^jz^k{\bar{z}}^l+O(|z|^5).
\]
Interestingly, these are the only coordinates which satisfy above vanishing conditions.\\
\begin{prop}
For any $m \in M$ there is unique holomorphic coordinate map $z = \phi_m$ up to affine transformation, for which there exist a \k potential $K_m(z)$ on $M$ such that all the $(0, l), \ (1, l), \ (l, 1), \ (l, 0)$ terms in the Taylor expansion of $K_m(z)$ vanish except for $(1,1)$ term which equal to $|z|^2$. These are exactly the K-coordinates.
\end{prop}
\begin{re}
It was in this form, K-coordinate first appeared in mathematics literature [B] and recently rediscovered in physics literature.\\
\end{re}
In above arguments, we did not worry about convergence of various power series, which in principle could cause problems. Now we would like to make some comments to address these problems.\\\\
Firstly, if the \k metric is actually real analytic or weakly real analytic (as defined later), which is usually true for many \k metrics that come up naturally, then those power series will converge in a neighborhood. In general, at least in principle, Those power series may not converge. But in most cases, we only need an approximate version of the canonical coordinates, namely coordinates that is ``canonical'' up to certain order, for which convergence is not a problem. Let the order goes to infinity, we will get a sequence of more and more accurate approximate ``canonical'' coordinates. These sequences of coordinates usually can do the same job as the canonical coordinates. For instance, this will be the case, when we later deal with the $C^\infty$ convergence of the Bergmann metrics.\\
\begin{de}
Let $f(x)$ be a smooth function defined in a neighborhood of zero, if Taylor series of $f(x)$ at $x=0$ converge in a neighborhood of $0$, then $f(x)$ is called weakly real analytic around $0$. $f(x)$ is called weakly real analytic function, if it is weakly real analytic at every point.
\end{de}
Clearly, real analytic functions are weakly real analytic. Following example is weakly real analytic at $0$, but not real analytic.\\\\
{\bf Example:} 
\[
f(x) = e^{-\frac{1}{x^2}} \ \ \ {\rm at} \ \ x=0
\]\\
Of course, it is an interesting question to see what kind of smooth functions are weakly real analytic. In the following, we will construct a smooth function, which is not weakly real analytic.\\\\
{\bf Example:}
\[
f(x) = \sum_{n=0}^\infty \frac{1}{(ib_n)^{n-1}} e^{i b_n x}
\]
Where
\[
b_n > {\rm Max}\left(2\sum_{m=0}^{n-1} (b_m)^{n-m+1}, (2n)!\right)
\]
With this assumption for $b_n$, we have
\[
f^{(n)}(x) = \sum_{m=0}^\infty (ib_m)^{n-m+1} e^{i b_m x}
\]
\begin{eqnarray*}
\left|f^{(n)}(x)\right| &\leq &\left|\sum_{m=0}^{n} (ib_m)^{n-m+1} e^{i b_m x}\right| + \left|\sum_{m=0}^\infty \frac{1}{(ib_{n+m +1})^m} e^{i b_{n+m+1} x}\right|\\
 &\leq &\left|\sum_{m=0}^{n} (ib_m)^{n-m+1} e^{i b_m x}\right| + 2
\end{eqnarray*}
Therefore $f(x)$ is smooth. On the other hand
\begin{eqnarray*}
\left|f^{(n)}(x)\right| &\geq &b_n - \left|\sum_{m=0}^{n-1} (ib_m)^{n-m+1} e^{i b_m x}\right| - \left|\sum_{m=0}^\infty \frac{1}{(ib_{n+m +1})^m} e^{i b_{n+m+1} x}\right|\\
&\geq &b_n - \sum_{m=0}^{n-1} (b_m)^{n-m+1} - \sum_{m=n+1}^\infty \frac{1}{((2m)!)^{m-n-1}}\\
&\geq & n^n n!
\end{eqnarray*}
Hence
\[
\left|\frac{f^{(n)}(x)}{n!}\right| \geq n^n
\]
Since 
\[
\sum_{n=0}^\infty n^n x^n
\]
is not convergent any where except at $x=0$, $f(x)$ is not weakly real analytic at any point, but it is smooth on the whole real line.\\\\ 
Now we assume that $\omega$ is an integer class, then there is an ample Hermitian line bundle $(L, h)$ over $(M, \omega)$ such that
\[
\omega = -\partial \bar{\partial} log h
\]
There exists a unique section $e_m(z)$ of $L|_{U_m}$ up to circle action such that $|e_m(0)|_h = 1$ and
\[
h(e_m(z), e_m(z)) = exp(-K_m)
\]
It is very interesting to look at the \ke metrics. When $g$ is \ke,
\[
Ric_g = - \partial \bar{\partial} log (det g) = \lambda g
\]
Since $(l, 0), (0, l)$ terms of g vanish, $(l, 1), (1, l)$ terms of $det g$ vanish. Although, presumably only $(l, 0), (0,l)$ terms of g vanish. Therefore, the natural local holomorphic section of canonical bundle of $M$ induced by determinant of the natural flat coordinate is exactly the unique section we discussed above for canonical bundle of $M$. The vanishing of $(l, 1), (1, l)$ terms of $det g$ merely reflect the fact that any trace with respect to $g$ of any $(l, 2), (2, l)$ terms of $K_m$ is zero.\\
\begin{re}
One thing we should keep in mind is that no positive definiteness of $\omega$ is required, so we may get affine coordinates with respect to any closed non-degenerate real analytic $(1,1)$ form, which is usually called pseudo-\k metric.\\
\end{re}
\begin{re}
Clearly all the $(l,k)$ components for $l,k > 2$ of $K_m$ are independent, they describe the local freedoms of \k manifolds. While for $(1,1)$ form in $\bf C^n$ the local freedom are characterized by $(l,k)$ components for $l,k \geq 2$.\\
\end{re}
\se{Tian's peak sections and Bergmann metrics}
The appearance of K-coordinate reminded me of an old paper of Tian's on Bergmann metric. In that paper, for a projective algebraic \k manifold $(M, g)$ with the associated ample Hermitian line bundle $(L, h)$ on $(M, g)$, he considered $L^m$ embedding $\varphi_m: M \longrightarrow {\bf CP^{N_m}}$, which is induced by an orthonormal basis $\{ S^m_0, \cdots , S^m_{N_m} \}$ of the space $H^0 (M, L^m)$ of all holomorphic global sections of $L^m$. Here the inner product on $H^0 (M, L^m)$ is the natural one induced by the \k metric $g$ and the Hermitian metric $h^m$ on $L^m$, i.e., $<S^m_\alpha, S^m_\beta> = \int_M h^m(S^m_\alpha, S^m_\beta) dV_g$. Let $g_{FS}$ be the standard Fubini-Study metric on ${\bf CP^{N_m}}$, i.e., $\omega_{g_{FS}} = \frac{\sqrt{-1}}{2}\partial \bar{\partial} log(\sum^{N_m}_{i=0} |w_i|^2)$. The pullback of the $\frac{1}{m}-$multiple of $g_{FS}$ on ${\bf CP^{N_m}}$, $g_m = \frac{1}{m}\varphi_m^* g_{FS}$ on $M$, is in the same \k class as $g$. It is called the Bergmann metric with respect to $L^m$.\\\\
It is important to notice that the Bergmann metric $g_m$ does not depend on the special choice of the orthonormal basis $\{ S^m_0, \cdots , S^m_{N_m} \}$. The main result in that paper is:\\\\
\begin{th}[Tian]
\[
\| g_m - g \|_{C^2} = O(\frac{1}{\sqrt{m}} )
\]
i.e. induced Bergmann metrics are $C^2$ convergent to the original metric.
\end{th}
The main point is to construct a sequence of peak sections of $L^m$ around a point using H\"{o}rmander $L^2$-estimate of $\bar{\partial}$-operator. Tian used a coordinate which is a lower order approximation of K-coordinate. The idea is in the local K-coordinate the canonical section $e_m(z)$ of $L$ has norm
\[
\|e_m(z)\| = \exp(-|z|^2 + \frac{1}{4}R_{i\bar{j}k\bar{l}}z^i z^{\bar{j}} z^k z^{\bar{l}}+O(|z|^5))
\]
\[
\|e^m_m(z)\| = \exp(-m(|z|^2 +O(|z|^4)))
\]
which is roughly the rapid degenerating Gauss function. When $m$ is large, it is very easy to smooth out and use H\"{o}rmander $L^2-$estimate of $\bar{\partial}-$operator to create a global holomorphic section that roughly looks like $e^m_m(z)$, which have the property that it is like a peak at $m$ and almost zero away from $m$.\\\\
One can actually create a sequence of sections with similar nature, such that for $P = (p_1, p_2, \cdots p_n)$
\[
e_P \sim z^{p_1}_1 z^{p_2}_2 \cdots z^{p_n}_n e^m_m
\]
Tian then went on using this sections to show his theorem. Due to the complexities that arise when order get higher, he only showed convergence for $C^2$, although he believed higher order estimate should be possible with a more neat method.\\\\
We will analyze the situation more closely and explore the $C^\infty$ convergence. First, let's fix some notations and summarize some known results.\\\\
For the notations, we will use the following convention, here letter ``P'' is just for illustration. We will also use other letters with similar convention.\\
For $P=(p_1, p_2, \cdots, p_n) \in {\bf Z^n}_+$, let
\[
p= \sum^n_{i=1} p_i, \ \ P!= \prod^n_{i=1} p_i !
\]
For $z=(z^1, z^2, \cdots, z^n) \in {\bf C^n}$, let
\[
dz= dz^1\wedge dz^2 \wedge \cdots \wedge dz^n,\ \ d\bar{z} = d\bar{z}^1\wedge d\bar{z}^2 \wedge \cdots \wedge d\bar{z}^n,\ \ dz\wedge d\bar{z}= \bigwedge^n_{i=1}dz^i\wedge d\bar{z} ^i
\]
\[
z^P = (z^1)^{p_1} (z^2)^{p_2}\cdots(z^n)^{p_n},\ \ \ \partial_P=(\frac{\partial}{\partial z^1})^{p_1}(\frac{\partial}{\partial z^2})^{p_2} \cdots (\frac{\partial}{\partial z^n})^{p_n}
\]
The construction of peak section in Tian's paper is well done, we really have no more to say on this. We simply state it as follows.
\begin{lm}[Tian]
For an n-tuple of integers $P=(p_1, p_2, \cdots, p_n) \in {\bf Z^n_+}$ and an integer $p'> p = p_1 + p_2 + \cdots + p_n$, there exists an $m_0>0$ such that for $m>m_0$, there is a holomorphic global section $S_P$ in $H^0 (M, L^m)$, satisfying
\[
\int_M\|S_P\|^2_{h^m} d V_g = 1,
\]
\[
\int_{M\backslash \{ \rho(z)\leq\frac{log m}{\sqrt{m}}\}} \|S_P\|^2_{h^m} d V_g = O(\frac{1}{m^{2p'}}),
\]
and locally at $x_0$,
\[
S_P(z) = \lambda_P(z^P + O(|z|^{2p'}))e_m^m(1+O(\frac{1}{m^{2p'}})),
\]
where $\|\cdot\|_{h^m}$ is the norm on $L^m$ given by $h^m$, and $O(1/m^{2p'})$ denotes a quantity dominated by $C/m^{2p'}$with the constant $C$ depending only on $p'$ and the geometry of $M$, moreover
\[
\lambda_P^{-2}= \int_{\rho(z)\leq\frac{log m}{\sqrt{m}}}|z^P|^2 a^m dV_g,
\]
where $dV_g = det(g_{i\bar{j}})(\sqrt{-1}/(2\pi))^n dz\wedge d\bar{z}$ is the volume form.
\end{lm}
For the estimates of the inner products between peak sections, we need the following straightforward generalization of Tian's.
\begin{lm}
\label{ca}
Let $S_P(z) = (1+O(\frac{1}{m^{2p'}}))(z^P + O(|z|^{2p'}))e_L^m$ be the section constructed as above. $T$ be another section of $L^m$.\\
(i) If $z^P$ is not in $T$'s Taylor expansion, then
\[
<S_P, T>_{h^m} = O(\frac{1}{m}) \|S_P\|_{h^m}\|T\|_{h^m}.
\]
(ii) If $T$ contains no terms $z^Q$, such that $q < p+d$ ($d \geq 1$) then
\[
<S_P, T>_{h^m} = O\left(\frac{1}{m^{1+d/2}}\right) \|S_P\|_{h^m}\|T\|_{h^m}.
\]
\end{lm}
{\bf Proof:}
We will only show $(ii)$, $(i)$ is simpler.\\\\
Assume that
\[
e^{m(|z|^2 - K(z, \bar{z}))} det g = \sum_{S, T} a_{S\bar{T}} z^S\bar{z}^T.
\]
Then
\[
<S_P, T>_{h^m}=\int_{|z|\leq \frac{log m}{\sqrt{m}}} z^P\overline{f_T}(z)e^{-m|z|^2}(\sum_{S, Q} a_{S\bar{Q}} z^S\bar{z}^Q) dzd\bar{z} + O(\frac{1}{m^{p'}}). 
\]
\[
=\int_{|z|\leq \frac{log m}{\sqrt{m}}} z^P\overline{f_T}(z)e^{-m|z|^2}(\sum_{s-q<d} a_{S\bar{Q}} z^S\bar{z}^Q) dz d\bar{z}
\]
\[
+ \int_{|z|\leq \frac{log m}{\sqrt{m}}} z^P\overline{f_T}(z) (O(m|z|^{d + 4} + |z|^{d+2})) h^m(z) dV_g
\]
\[
\leq C \left[\int_{|z|\leq \frac{log m}{\sqrt{m}}} |z^P|^2 (m^2 |z|^{2d + 8} + |z|^{2d+4}) h^m(z) dV_g\right]^{\frac{1}{2}} \| T \|_{h^m}
\]
\[
= O\left(\frac{1}{m^{\frac{d+2}{2}}}\right)\| S_P \|_{h^m}\| T \|_{h^m}
\]
{\bf Q.E.D.}\\\\
Now we will compute inner products of $S_P$. Here we make a modification, replacing Tian's factor $(1-|z|^2)^m$ by $e^{-m|z|^2}$ which seems more natural and easier to compute.\\\\
\begin{lm}
\[
(\frac{\sqrt{-1}}{2})^n \int_{\rho(z)\leq\frac{log m}{\sqrt{m}}}|z^P|^2 e^{-m|z|^2} dz\wedge d\bar{z}=(\frac{\pi}{m})^n \frac{P!}{m^p} + O(\frac{1}{m^{2p'}}) \ \ {\rm for}\ m \ {\rm large}
\]
\end{lm}
{\bf Proof:}
It is easy to see that it will not make a difference, if we extend domain of integration to ${\bf C^n}$. Then the integral is of separated variables and result is a corollary of the fact:
\[
\int^\infty_0 r^{2p} e^{-mr^2}rdr = \frac{1}{2}\frac{p!}{m^p}
\]
{\bf Q.E.D.}\\\\
Now recall that
\[
e^{m(|z|^2 - K(z, \bar{z}))} det g = \sum_{S, T} a_{S\bar{T}} z^S\bar{z}^T.
\]
Then we have the following lemma.
\begin{lm}
\label{cb}
\[
<S_P, S_Q>_{h^m} = (\frac{\pi}{m})^n \sum_{S+P=Q+T}a_{S\bar{T}} \frac{(S+P)!}{m^{s+p}}.
\]
Especially
\[
\|S_P\|^2_{h^m} = (\frac{\pi}{m})^n \sum_{S}a_{S\bar{S}} \frac{(S+P)!}{m^{s+p}}
\]
with the leading term 
\[
(\frac{\pi}{m})^n\frac{P!}{m^{p}}.
\]\\
\end{lm}
Notice:
\[
\omega_m =  \frac{\sqrt{-1}}{2\pi}\partial \bar{\partial}\frac{1}{m} log(\sum^{N_m}_{i=0} |S^m_i|^2)
\]
\[
\omega_g =  \frac{\sqrt{-1}}{2\pi}\partial \bar{\partial}K
\]
here
\[
K = |z|^2 - \frac{1}{4} R_{i\bar{j}k\bar{l}}z^i\bar{z}^j z^k\bar{z}^l + O(|z|^5)
\]
For us it is important to notice that $\omega_m$ converge to $\omega_g$ in $C^\infty$ is equivalent to that the \k potentials $K_m = \frac{1}{m} log(\sum^{N_m}_{i=0} |S^m_i|^2)$ converge to $K$ in $C^\infty$. This is what we will try to show.\\\\
First put index $P$ in an order, let $|P|$ denote its order. The idea is if we care about the $z^P$ term, we want to take a $D \gg P$ and also $m \gg d, p' \gg d$ such that, we can just use $S_P, \  P\leq D$ to do the computation and be able to ignore terms $>D$. This need a very careful arrangement.\\\\
First construct a filtration of $H^0 (M, L^m)$, let 
\[
V_P = \{S \in H^0 (M, L^m): \partial_Q S|_m = 0 \ {\rm for}\ Q<P\}
\]
If $Q>P$ is ``smallest'', then $dim(V_Q/V_P) \leq 1$.\\\\
We can then associate an orthonormal basis of $H^0 (M, L^m)$ with respect to this filtration, denote by $\{ S^m_0, \cdots , S^m_{N_m} \}$. We replace this by $\{ T_0, T_{e_1}, \cdots T_D,$ $S^m_{|D|+1}, \cdots , S^m_{N_m} \}$ where $T_P=S_P/\|S_P\|_{h^m}$.\\\\
We intend to make them orthogonal again, but in a special way. We first modify $T_P$ for $P\leq D$ to make them orthogonal to $S^m_i$ for $i>|D|$. By which we get
\[
\tilde{T}_P = T_P + \sum^{N_m}_{j=|D|+1} \beta_{Pj} S^m_j
\]
Use lemma \ref{ca} for $S_P$ and $\sum^{N_m}_{j=|D|+1} \beta_{Pj} S^m_j$, we get
\[
\sum^{N_m}_{j=|D|+1} |\beta_{Pj}|^2 = O(m^{-(d+3-p)})
\]
Now we try to rectify $\{\tilde{T}_P\}$. Our aim is to rule out the influence of $\{S^m_j\}^{N_m}_{j=|D|+1}$. We will find out, after all, if we care $z^P$ term, the influence of $\{S^m_j\}$ only reflect on the coefficients of $\beta_{QR}$, $Q, R < P$. Since all the other terms vanish at $x_0$. Now we try to trace the influence of $\{S^m_j\}$ on $\beta_{QR}$. Rectifying $\{\tilde{T}_P\}$ requires operations similar to
\[
\tilde{T}_P - \frac{<\tilde{T}_P, \tilde{T}_Q>}{|\tilde{T}_P||\tilde{T}_Q|}\tilde{T}_Q
\]
where $Q>P$. The major terms coming from $\{S^m_j\}$ is
\[
- <T_P, \sum^{N_m}_{j=|D|+1} \beta_{Qj} S^m_j> - <T_Q, \sum^{N_m}_{j=|D|+1} \beta_{Pj} S^m_j> 
\]
\[
+ \frac{1}{2} <T_P, T_Q><T_Q, \sum^{N_m}_{j=|D|+1} \beta_{Qj} S^m_j>
\]
which is of order $O( m^{\frac{p+q}{2}-(3+d)})$. A careful analysis will show that other terms are of lower order and can be ignored.\\\\
Assume that after rectification, we get
\[
\hat{T_P} = T_P + \sum_{P<Q\leq D} \beta_{PQ} T_Q + \sum^{N_m}_{j=|D|+1} \beta_{Pj} S^m_j
\]
For simplicity, here we still use the same notation for $\beta_{Pj}$, although it have been changed. It is easy to see its order is unchanged. Then we have 
\[
\beta_{PQ} = \gamma_{PQ} + O( m^{\frac{p+q}{2}-(3+d)})
\]
Here $\gamma_{PQ}$ comes from pure interaction of $\{ T_P \}$.\\\\
We are interested in 
\[
\frac{1}{m} log\left(\sum^{N_m}_{i=0} |S^m_i|^2\right)=\frac{1}{m} log\left(\sum_{0\leq Q \leq D} \frac{|\hat{T_Q}|^2}{\|\hat{T_Q}\|^2_{h^m}} + \sum^{N_m}_{i=|D|+1} |S^m_i|^2\right)
\]
When we worry about $z^P$ term, only ${|\hat{T_Q}|^2}/{\|\hat{T_Q}\|^2_{h^m}}$ terms for $Q \leq P$ are involved. A careful analysis will reveal that ${\|\hat{T_Q}\|^2_{h^m}}$ terms have error terms controlled by $O(m^{p-(4+d)})$, and the ${|\hat{T_Q}|^2}$ terms have error terms controlled by $O(m^{p-(3.5+d)})$. They are all negligible, since $d \gg p$.\\\\
Above argument ensure that we only need to work with the subspace of $H_0(M, L^m)$ generated by $\{S_P\}_{p\leq d}$. Another very crucial observation is that although orthonormal basis is good for abstract arguments, it is not suitable for computation. We will use the more natural but non-orthogonal basis $\{S_P\}$ to do the computation.\\\\
Let $\{T_i^m\}$ be a not necessarily orthonormal basis of $H^0(M, L^m)$, and $\{S_i^m\}$ be the original orthonormal basis. Then
\[
|S^m_0|^2 + |S^m_1|^2 + \cdots |S^m_{N_m}|^2 = \sum^{N_m}_{i,j=0}G^{i \bar{j}} T^m_i \overline{T^m_j}
\]
Here $G^{i \bar{j}}G_{k \bar{j}}= \delta_{ik}$, $G_{i \bar{j}}= <T^m_i, T^m_j>_{h^m}$.\\\\
So we only need to compute
\[
\frac{1}{m} log\left(\sum_{P, Q} G^{P\bar{Q}} S_P \overline{S_Q}\right)
\]
Here $G^{P \bar{Q}}G_{R \bar{Q}}= \delta_{PR}$, $G_{P \bar{Q}}= <S_P, S_Q>_{h^m}$.\\\\
Let \k potential be
\[
K(z, \bar{z}) = |z|^2 + \sum_{P,Q} c_{P \bar{Q}} z^P \bar{z}^Q
\]
What we expect is  that all $m^r$ terms for $\frac{1}{m} log(\sum_{P, Q} G^{P\bar{Q}} S_P \overline{S_Q})$ of $r>0$ vanish and the terms that do not involve $m$ should equal to $K(z, \bar{z})$.\\\\
Notice that $\{c_{P \bar{Q}}\}$ are independent, so we may trace them independently. Now we have reduced our problem into a purely combinatoric one. In principle, the problem is already solved, if one is optimistic. But it turns out the corresponding combinatoric problem is extremely complicated. Although we are able to sort out almost all combinatorics, the complete solution of the problem needs to go back to some geometrical arguments.\\\\
\se{Combinatorics Involved}
In this section we will demonstrate the combinatorics involved, and try to solve it to a great extent in next section.\\\\
As we mentioned before that $\{c_{P \bar{Q}}\}$ are independent, so without loss of generality, we only need to consider finitely many terms. Let \k potential be\\
\[
K(z, \bar{z}) = |z|^2 + \sum_{i=1}^r c_{P^i \bar{Q}^i} z^{P^i} \bar{z}^{Q^i}
\]
then\\ 
\[
g_{i\bar{j}} = \delta_{ij} + \sum_{k=1}^r p^k_{,i}q^k_{,j}c_{P^k \bar{Q}^k} z^{P^k-e_i} \bar{z}^{Q^k-e_j} 
\]
and\\
\begin{eqnarray*}
\det (g) &= &1 + \sum_{i=1}^r \sum_{j=1}^n p^i_{,j}q^i_{,j}c_{P^i \bar{Q}^i} z^{P^i-e_j} \bar{z}^{Q^i-e_j} + \cdots\\
&=&1 + \sum_{i=1}^r c_{P^i \bar{Q}^i} z^{P^i} \bar{z}^{Q^i}\sum_{j=1}^n \left( \frac{p^i_{,j}}{z^j}\right) \left( \frac{q^i_{,j}}{\bar{z}^j}\right)\\
&&+ \sum_{I=(i_1,i_2)\subset (1,2,\cdots,r)} c_{P^{i_1} \bar{Q}^{i_1}}c_{P^{i_2} \bar{Q}^{i_2}} z^{P^{i_1}+P^{i_2}}\bar{z}^{Q^{i_1}+Q^{i_2}} \\
&&\sum_{J=(j_1,j_2)\subset (1,2,\cdots,n)} \frac{\det(P^I_{,J})}{z^J} \frac{\det(Q^I_{,J})}{\bar{z}^J} + \cdots\\
&=& \sum_{I\subset (1,2,\cdots,r)} \left( \prod_{i\in I} c_{P^i \bar{Q}^i} z^{P^i} \bar{z}^{Q^i}\right)\sum_{J\subset (1,2,\cdots,n)} \frac{\det(P^I_{,J})}{z^J} \frac{\det(Q^I_{,J})}{\bar{z}^J}
\end{eqnarray*} 
Recall\\
\[
e^{m(|z|^2 - K(z, \bar{z}))} det g = \sum_{S, T} \tilde{a}_{S\bar{T}} z^S\bar{z}^T.
\]
Here we use $\tilde{a}_{S\bar{T}}$ instead of $a_{S\bar{T}}$ which was used. Let\\
\[
b_{P\bar{Q}} = - c_{P\bar{Q}}m^{1-\frac{p+q}{2}}, \ \ a_{P\bar{Q}}=\tilde{a}_{P\bar{Q}}m^{-\frac{p+q}{2}}
\]
Then\\
\[
\sum_{S, T} a_{S\bar{T}} (\sqrt{m}z)^S(\sqrt{m}\bar{z})^T=\sum_{S, T} \tilde{a}_{S\bar{T}} z^S\bar{z}^T =  e^{m(|z|^2 - K(z, \bar{z}))} \det g
\]
\begin{eqnarray*}
&=&e^{-m(\sum_{i=1}^r c_{P^i \bar{Q}^i} z^{P^i} \bar{z}^{Q^i})} \left(1 + \sum_{i=1}^r c_{P^i \bar{Q}^i}\sum_{j=1}^n p^i_{,j}q^i_{,j} z^{P^i-e_j} \bar{z}^{Q^i-e_j} + \cdots \right)\\
&=& \sum_{P, Q} \left(\sum_{\tiny{\begin{array}{c} L=(l_1,l_2, \cdots l_r)\\ \sum_{i=1}^r l_i(P^i, Q^i)=(P, Q)\end{array}}} \frac{1}{L!}\prod^r_{i=1} b_{P^i\bar{Q}^i}^{l_i}\right)(\sqrt{m}z)^P(\sqrt{m}\bar{z})^Q\\
&&\left(1 - \sum_{i=1}^r b_{P^i \bar{Q}^i}\sum_{j=1}^n p^i_{,j}q^i_{,j} (\sqrt{m}z)^{P^i-e_j} (\sqrt{m}\bar{z})^{Q^i-e_j} + \cdots \right)\\
&=& \sum_{P, Q} \left(\sum_{\tiny{\begin{array}{c} L=(l_1,l_2, \cdots l_r)\\ \sum_{i=1}^r l_i(P^i, Q^i)=(P, Q)\end{array}}} \frac{1}{L!}\prod^r_{i=1} b_{P^i\bar{Q}^i}^{l_i}\right)(\sqrt{m}z)^P(\sqrt{m}\bar{z})^Q\\
&&\sum_{I\subset (1,2,\cdots,r)} (-1)^{|I|}\left( \prod_{i\in I} b_{P^i \bar{Q}^i} (\sqrt{m}z)^{P^i} (\sqrt{m}\bar{z})^{Q^i}\right)\sum_{J\subset (1,2,\cdots,n)} \frac{\det(P^I_{,J})}{(\sqrt{m}z)^J } \frac{\det(Q^I_{,J})}{(\sqrt{m}\bar{z})^J}
\end{eqnarray*}
Let $S_P = z^Pe_L^m$ be the sections we constructed. Then according to lemma \ref{cb} we have\\
\[
<S_P, S_Q>_{h^m} = (\frac{\pi}{m})^n \sum_{S+P=Q+T}\tilde{a}_{S\bar{T}} \frac{(S+P)!}{m^{s+p}}.
\]
We define\\
\[
G_{S\bar{T}} = \sum_{S+P=Q+T}\tilde{a}_{P\bar{Q}} \frac{(S+P)!}{m^{s+p}} = \frac{S!}{m^{\frac{s+t}{2}}}\left(\sum_{S+P=Q+T}  a_{P\bar{Q}} \frac{(S+P)!}{S!}\right)
\]
Let\\ 
\[
\Lambda_{S\bar{T}} = \frac{S!}{m^s} \delta_{ST} \ \ \Lambda^{S\bar{T}} = \frac{m^s}{S!} \delta_{ST}
\]
Then
\[
G_{S\bar{T}} = \Lambda_{S\bar{T}} + G_{S\bar{T}}'
\]\\
To compute inverse of matrix $G = (G_{S\bar{T}}) = \Lambda + G'$, we use the following formula.\\
\[
G^{-1} = \Lambda^{-1} - \Lambda^{-1}G'\Lambda^{-1} + (\Lambda^{-1}G')^{2}\Lambda^{-1} + \cdots
\]
Use this equality we get
\begin{eqnarray*}
G^{T\bar{S}}& =& \Lambda^{S\bar{T}} + \frac{m^{\frac{s+t}{2}}}{S!}\sum_{u=1}^\infty (-1)^u \sum_{\tiny{\sum_{j=1}^u Q_j - P_j = S -T}}\\
&&\left(\prod_{j=1}^u a_{P_j \bar{Q_j}}\right)\left(\prod_{j=1}^u\frac{(S+ \sum_{i=1}^j (P_i - Q_i) + Q_j)!}{(S+ \sum_{i=1}^j (P_i - Q_i))!}\right)
\end{eqnarray*}\\
Here $Q_j, P_j \in \bf{Z} _+^n$.\\
To write in terms of $b$, we recall,
\[
\sum_{S, T} a_{S\bar{T}} (\sqrt{m}z)^S(\sqrt{m}\bar{z})^T=\sum_{S, T} \tilde{a}_{S\bar{T}} z^S\bar{z}^T
\]
\begin{eqnarray*}
&=& \sum_{P, Q} \left(\sum_{\tiny{\begin{array}{c} L=(l_1,l_2, \cdots l_r)\\ \sum_{i=1}^r l_i(P^i, Q^i)=(P, Q)\end{array}}} \frac{1}{L!}\prod^r_{i=1} b_{P^i\bar{Q}^i}^{l_i}\right)(\sqrt{m}z)^P(\sqrt{m}\bar{z})^Q\\
&&\left(1 - \sum_{i=1}^r b_{P^i \bar{Q}^i}\sum_{j=1}^n p^i_{,j}q^i_{,j} (\sqrt{m}z)^{P^i-e_j} (\sqrt{m}\bar{z})^{Q^i-e_j} + \cdots \right)\\
&=& \sum_{P, Q} \left(\sum_{\tiny{\begin{array}{c} L=(l_1,l_2, \cdots l_r)\\ \sum_{i=1}^r l_i(P^i, Q^i)=(P, Q)\end{array}}} \frac{1}{L!}\prod^r_{i=1} b_{P^i\bar{Q}^i}^{l_i}\right)(\sqrt{m}z)^P(\sqrt{m}\bar{z})^Q\\
&&- \sum_{P, Q} \sum_{\tiny{\begin{array}{c} L=(l_1,l_2, \cdots l_r)\\ \sum_{i=1}^r l_i(P^i, Q^i)=(P, Q)\end{array}}} \left( \prod^r_{i=1} b_{P^i\bar{Q}^i}^{l_i}\right) \\
&&\sum_{i=1}^r \sum_{j=1}^n p^i_{,j}q^i_{,j}\frac{1}{(L-e_i)!} (\sqrt{m}z)^{P-e_j} (\sqrt{m}\bar{z})^{Q-e_j}\cdots\\
&=&\sum_{P, Q} \sum_{\tiny{\begin{array}{c} L=(l_1,l_2, \cdots l_r)\\ \sum_{i=1}^r l_i(P^i, Q^i)=(P, Q)\end{array}}} \left( \prod^r_{i=1} b_{P^i\bar{Q}^i}^{l_i}\right) \\
&&\sum_{I\subset (1,2, \cdots, r)} \sum_{J\subset (1,2, \cdots, n)} \det(P^I_{,J})\det(Q^I_{,J})\frac{(-1)^{|I|}}{(L-e_I)!} (\sqrt{m}z)^{P-e_J} (\sqrt{m}\bar{z})^{Q-e_J} 
\end{eqnarray*}\\
which gives
\begin{eqnarray*}
a_{P\bar{Q}}&=&\left(\sum_{\tiny{\begin{array}{c} L=(l_1,l_2, \cdots l_r)\\ \sum_{i=1}^r l_i(P^i, Q^i)=(P, Q)\end{array}}} \frac{1}{L!}\prod^r_{i=1} b_{P^i\bar{Q}^i}^{l_i}\right)\\
&&-\sum_{j=1}^n\sum_{\tiny{\begin{array}{c} L=(l_1,l_2, \cdots l_r)\\ \sum_{i=1}^r l_i(P^i, Q^i)=(P+e_j, Q+e_j)\end{array}}} \left( \prod^r_{i=1} b_{P^i\bar{Q}^i}^{l_i}\right) \sum_{i=1}^r p^i_{,j}q^i_{,j}\frac{1}{(L-e_i)!} + \cdots\\
&=&\sum_{J\subset (1,2, \cdots, n)}\sum_{\tiny{\begin{array}{c} L=(l_1,l_2, \cdots l_r)\\ \sum_{i=1}^r l_i(P^i, Q^i)=(P+e_J, Q+e_J)\end{array}}} \left( \prod^r_{i=1} b_{P^i\bar{Q}^i}^{l_i}\right)\\
&&\sum_{I\subset (1,2, \cdots, r)} \det(P^I_{,J})\det(Q^I_{,J})\frac{(-1)^{|I|}}{(L-e_I)!}
\end{eqnarray*}\\
Then 
\begin{eqnarray*}
G^{T\bar{S}}& =& \Lambda^{S\bar{T}} + \frac{m^{\frac{s+t}{2}}}{S!}\sum_{u=1}^\infty (-1)^u \sum_{\tiny{\begin{array}{c} L=(l_1,l_2, \cdots l_r)\\ \sum_{j=1}^u L^j = L \\ \sum_{i=1}^r l_i(P^i - Q^i)=T-S\end{array}}} \frac{1}{\prod^u_{j=1} {(L^j)!}}\\
&&\left( \prod^r_{i=1} b_{P^i\bar{Q}^i}^{l_i}\right) \left(\prod_{j=1}^u\frac{(S+ \sum_{i=1}^j (P_i - Q_i) + Q_j)!}{(S+ \sum_{i=1}^j (P_i - Q_i))!}\right)\\
&&\left(1-\sum_{j=1}^u \sum_{i=1}^r \sum_{k=1}^n p^i_{,k}q^i_{,k}l_i^j \frac{1}{s_{,k} + \sum_{i=1}^j(p_{i,k} - q_{i,k}) + q_{j,k}} + \cdots \right)\\\\
&=& \Lambda^{S\bar{T}} + \frac{m^{\frac{s+t}{2}}}{S!}\sum_{u=1}^\infty (-1)^u \sum_{\tiny{\begin{array}{c} L=(l_1,l_2, \cdots l_r)\\ \sum_{j=1}^u L^j = L \\ \sum_{i=1}^r l_i(P^i - Q^i)=T-S\end{array}}} \frac{1}{\prod^u_{j=1} {(L^j)!}}\left( \prod^r_{i=1} b_{P^i\bar{Q}^i}^{l_i}\right) \\
&&\left(\prod_{j=1}^u\frac{(S+ \sum_{i=1}^j (P_i - Q_i) + Q_j)!}{(S+ \sum_{i=1}^j (P_i - Q_i))!}\right)\left(\prod_{j=1}^u \sum_{\tiny{\begin{array}{c}{I_j\subset (1,2, \cdots, r)}\\ {K_j\subset (1,2, \cdots, n)}\end{array}}} \right.\\
&& \left. \det(P^{I_j}_{,K_j})\det(Q^{I_j}_{,K_j}) (L^j)^{e_{I_j}} \frac{(-1)^{|K_j|}}{\prod_{k\in K_j}(s_{,k} + \sum_{i=1}^j(p_{i,k} - q_{i,k}) + q_{j,k})}  \right)
\end{eqnarray*}\\
Here $e_I = \sum_{i\in I} e_i$, $L^j = (l^j_1, l^j_1, \cdots l^j_r)$, and $P_j = \sum^r_{i=1} l^j_i P^i$, $Q_j = \sum^r_{i=1} l^j_i Q^i$. Let\\
\[
I_L = \sum_{u=1}^\infty (-1)^u \sum_{\tiny{\sum_{j=1}^u L^j = L }} \frac{1}{\prod^u_{j=1} {(L^j)!}}
\]
Here $L = (l_1, l_1, \cdots l_r)$. Notice that\\
\[
\sum_L I_L T^L = \sum_{u=1}^\infty (-1)^u (e^t - 1)^u = e^{-t} = \sum_L \frac{(-1)^l}{L!} T^L
\]
here $T = (t_1, t_1, \cdots t_r), \ t= (t_1+ t_1+ \cdots +t_r)$. Therefore\\
\begin{equation}
I_L = \sum_{u=1}^\infty (-1)^u \sum_{\tiny{\sum_{j=1}^u L^j = L }} \frac{1}{\prod^u_{j=1} {(L^j)!}} = \frac{(-1)^l}{L!}
\label{da}
\end{equation}
With this equality in mind, we can write $G^{T\bar{S}}$ as
\[
G^{T\bar{S}} = \Lambda^{S\bar{T}} + \frac{m^{\frac{s+t}{2}}}{S!} \sum_{\tiny{\begin{array}{c} L=(l_1,l_2, \cdots l_r) \\ \sum_{i=1}^r l_i(P^i - Q^i)=T-S\end{array}}}\frac{\prod^r_{i=1} (-b_{P^i\bar{Q}^i})^{l_i}}{L!}
\]
\[
\left(\frac{S!}{(S- \sum_{i=1}^r l_iQ^i)!} + \sum_{j=1}^n D_{1,j}(L, P, Q)\frac{S!}{(S- \sum_{i=1}^r l_iQ^i + e_j)!}+ \cdots \right)
\]
Use this expression to compute the \k potential, we get
\[
K_m = \frac{1}{m} log \left(1 + \sum_{S, T}\frac{G^{T\bar{S}}}{G^{0\bar{0}}}\bar{z}^S z^T\right)
\]
\[
= \frac{1}{m} \sum_{k=1}^\infty \frac{(-1)^{k-1}}{k} \sum_{S, T} \bar{z}^S z^T \sum_{\tiny{\begin{array}{c} \sum_{i=1}^k S_i = S \\ \sum_{i=1}^k T_i = T\end{array}}} \prod_{i=1}^k \frac{G^{T_i\bar{S}_i}}{G^{0\bar{0}}}
\]
\[
= \frac{1}{m} \sum_{k=1}^\infty \frac{(-1)^{k-1}}{k} \sum_{S, T} (\sqrt{m}z)^S(\sqrt{m}\bar{z})^T  \sum_{\tiny{\begin{array}{c} \sum_{i=1}^k S_i = S \\ \sum_{i=1}^k T_i = T\end{array}}} \frac{1}{\prod_{i=1}^k S_i!}
\]
\begin{equation}
\label{db}
\left[\sum_{\tiny{\begin{array}{c} L=(l_1,l_2, \cdots l_r) \\ \sum_{i=1}^r l_i(P^i - Q^i)=T-S\end{array}}}\prod^r_{i=1} (-b_{P^i\bar{Q}^i})^{l_i}\sum_{\tiny{\sum_{i=1}^k L^i = L }} \frac{1}{\prod^k_{i=1} {(L^i)!}}\right.
\end{equation}
\[
\left.\prod_{i=1}^k\left(\frac{S_i!}{(S_i- \sum_{j=1}^r l^i_jQ^j)!} + \sum_{v=1}^n D_{1,v}(L, P_v, Q_v)\frac{S_i!}{(S_i- \sum_{j=1}^r l_j^iQ^j + e_v)!}+ \cdots \right)\right] 
\]
\[
= \frac{1}{m}\sum_{S-T}\sum_{\tiny{\begin{array}{c} L=(l_1,l_2, \cdots l_r) \\ \sum_{i=1}^r l_i(P^i - Q^i)=T-S\end{array}}}\prod^r_{i=1} (-b_{P^i\bar{Q}^i})^{l_i}\sum_{u=1}^\infty \sum_{\tiny{\sum_{i=1}^u L^i = L }} \frac{1}{\prod^u_{i=1} {(L^i)!}}
\]
\[
\sum_S (\sqrt{m}z)^S(\sqrt{m}\bar{z})^T \sum_{k=u}^\infty \frac{(-1)^{k-1}}{k}\sum_{\tiny{ \sum_{i=1}^k S_i = S }} \frac{1}{\prod_{i=1}^k S_i!}
\]
\[
\left(\sum _{\{j_i\}_{i=1}^u \in \{1,2, \cdots k\}} \prod_{i=1}^u\left(\frac{S_{j_i}!}{(S_{j_i}- \sum_{k=1}^r l^i_k Q^k)!} + \cdots\right)\right)
\]
To this point, we have transformed the problem to a clear combinatoric one. What we need to do is to reduce the above expression via combinatoric techniques to $K + O(\frac{1}{m})$. We can either tackle the combinatoric problem head on, or try to get around it in some way. In the next section, we will see how much we can do with pure combinatoric. Then in the further section, we will present a rather simple geometric proof.\\\\
\se{Combinatoric Approach}
To proceed further, we need some more sophisticated combinatoric equalities. Recall from $(\ref{da})$, we have
\[
\sum_{u=1}^\infty (-1)^u \sum_{\tiny{\sum_{j=1}^u L^j = L }} \frac{1}{\prod^u_{j=1} {(L^j)!}} = \frac{(-1)^l}{L!}
\]
Use a similar method we can get the following.\\
\begin{lm}
For $L \in {\bf Z^n_+}$ we have
\[
\sum_{u=1}^\infty \frac{(-1)^{u-1}}{u} \sum_{\tiny{\sum_{j=1}^u L^j = L }} \frac{1}{\prod^u_{j=1} {(L^j)!}} = 
\left\{
\begin{array}{c}
0, \ \ l >1\\
1, \ \ l =1
\end{array}
\right.
\]\\
\end{lm}
{\bf Proof:} Let
\[
I_L= \sum_{u=1}^\infty \frac{(-1)^{u-1}}{u} \sum_{\tiny{\sum_{j=1}^u L^j = L }} \frac{1}{\prod^u_{j=1} {(L^j)!}}
\]\\
Notice that 
\[
\sum_L I_L T^L = \sum_{u=1}^\infty \frac{(-1)^{u-1}}{u} (e^t - 1)^u = t
\]\\
The equation in the lemma then follow.\\
{\bf Q.E.D.}\\\\
With the above two equalities in mind, we will show some generalizations of them. The method is similar, but need more care. The following elementary equality is used in the proofs.\\\\
\[
(1+x)^{-r} = \sum_u^\infty (-1)^{u+r-1} \frac{(u+r-1)!}{u!(r-1)!} x^u
\]\\
\begin{lm}
\label{ec}
For $L, A_1, A_2, \cdots A_r \in {\bf Z^n_+}$ and $A = \sum_{i=1}^r A_i$ we have
\[
\sum_{u=r}^\infty (-1)^u \sum_{\tiny{\sum_{j=1}^u L^j = L }} \frac{1}{\prod^u_{j=1} {(L^j)!}}\left(\sum^u_{\tiny{\begin{array}{c} j_k = 1\\ k=1,2, \cdots r \end{array}}} \prod^r_{k=1} \frac{L^{j_k}!}{(L^{j_k} - A_k)!}\right) 
\]
\[
= \left\{
\begin{array}{cc}
(-1)^{r+l-a}\frac{r!}{(L-A)!}, \ & \ L \geq A\\
0,& \ L < A
\end{array}
\right.
\]
\end{lm}
{\bf Proof:} Let
\[
g(L,u) = \sum_{\tiny{\sum_{j=1}^u L^j = L }} \frac{1}{\prod^u_{j=1} {(L^j)!}}
\]
Then
\[
\sum_L t^L \sum_{u=r}^\infty (-1)^u \sum_{\tiny{\sum_{j=1}^u L^j = L }} \frac{1}{\prod^u_{j=1} {(L^j)!}}\left(\sum^u_{\tiny{\begin{array}{c} j_k = 1\\ k=1,2, \cdots r \end{array}}} \prod^r_{k=1} \frac{L^{j_k}!}{(L^{j_k} - A_k)!}\right) 
\]
\begin{eqnarray*}
&=&\sum_L t^L \sum_{u=r}^\infty (-1)^u \frac{u!}{(u-r)!} \sum_{i=0}^r\left( \begin{array}{c}r\\i\end{array}\right) g(L-A, u-i)\\
&=&\sum_u^\infty (-1)^u \frac{u!}{(u-r)!} \sum_{i=0}^r \left( \begin{array}{c}r\\i\end{array}\right)\sum_L t^L g(L-A, u-i)\\
&=&\sum_u^\infty (-1)^u \frac{u!}{(u-r)!} t^A \sum_{i=0}^r \left( \begin{array}{c}r\\i\end{array}\right)(e^t -1)^{u-i}\\
&=&t^A e^{tr} \sum_u^\infty (-1)^u \frac{u!}{(u-r)!} (e^t - 1)^{u-r}\\
&=&(-1)^r r! e^{-t} t^A\\
\end{eqnarray*}
Compare the coefficients, we get the lemma.\\\\
{\bf Q.E.D.}\\\\
\begin{lm}
\label{ea}
For $L, A_1, A_2, \cdots A_r \in {\bf Z^n_+}$ and $A = \sum_{i=1}^r A_i$ we have\\\\
(i) when $r > 0$
\[
\sum_{u=r}^\infty \frac{(-1)^{u-1}}{u} \sum_{\tiny{\sum_{j=1}^u L^j = L }} \frac{1}{\prod^u_{j=1} {(L^j)!}}\left(\sum^u_{\tiny{\begin{array}{c} j_k = 1\\ k=1,2, \cdots r \end{array}}} \prod^r_{k=1} \frac{L^{j_k}!}{(L^{j_k} - A_k)!}\right) 
\]
\[
= \left\{
\begin{array}{cc}
(-1)^{r-1}(r-1)!, \ & \ L = A\\
0,& \ L \not= A
\end{array}
\right.
\]
(ii) when $r = 0$
\[
\sum_{u=1}^\infty \frac{(-1)^{u-1}}{u} \sum_{\tiny{\sum_{j=1}^u L^j = L }} \frac{1}{\prod^u_{j=1} {(L^j)!}}
= \left\{
\begin{array}{cc}
1, \ & \ l = 1\\
0, \ & \ l \not= 1
\end{array}
\right.
\]
\end{lm}
{\bf Proof:} Case $(ii)$ is simple, we will prove case $(i)$, so assume $r > 0$.
\[
\sum_L t^L \sum_{u=r}^\infty \frac{(-1)^{u-1}}{u} \sum_{\tiny{\sum_{j=1}^u L^j = L }} \frac{1}{\prod^u_{j=1} {(L^j)!}}\left(\sum^u_{\tiny{\begin{array}{c} j_k = 1\\ k=1,2, \cdots r \end{array}}} \prod^r_{k=1} \frac{L^{j_k}!}{(L^{j_k} - A_k)!}\right) 
\]
\begin{eqnarray*}
&=&\sum_L t^L \sum_{u=r}^\infty \frac{(-1)^{u-1}}{u} \frac{u!}{(u-r)!} \sum_{i=0}^r\left( \begin{array}{c}r\\i\end{array}\right) g(L-A, u-i)\\
&=&\sum_u^\infty \frac{(-1)^{u-1}}{u} \frac{u!}{(u-r)!} \sum_{i=0}^r \left( \begin{array}{c}r\\i\end{array}\right)\sum_L t^L g(L-A, u-i)\\
&=&\sum_u^\infty \frac{(-1)^{u-1}}{u} \frac{u!}{(u-r)!} t^A \sum_{i=0}^r \left( \begin{array}{c}r\\i\end{array}\right)(e^t -1)^{u-i}\\
&=&t^A e^{tr} \sum_u^\infty (-1)^{u-1} \frac{(u-1)!}{(u-r)!} (e^t - 1)^{u-r}\\
&=&(r-1)! t^A e^{tr} \sum_u^\infty (-1)^{u+r-1} \frac{(u+r-1)!}{u!(r-1)!} (e^t - 1)^u\\
&=&(-1)^{r-1} (r-1)! t^A e^{tr}(\frac{1}{e^t})^r\\
&=&(-1)^{r-1} (r-1)! t^A
\end{eqnarray*}
Compare the coefficients, we get the lemma.\\\\
{\bf Q.E.D.}\\\\
In our proof, we will also need following type equality. The idea is to express polynomials of $s$ in term of polynomials like: $\frac{s!}{(s-m)!}$ for $m \in {\bf Z_+}$.\\\\
\begin{lm}
\[
\frac{(s+p)!}{(s+p-q)!} = \sum_{i=0} \left(\begin{array}{c}p\\i\end{array}\right)\left(\begin{array}{c}q\\i\end{array}\right)\frac{i!s!}{(s-q+i)!}
\]
\end{lm}
{\bf Proof:} It is easy to show by induction.\\\\
{\bf Q.E.D.}\\\\
We have the following multivariable version of this lemma, notice that the multi-index $I \in {\bf Z^n_+}$.
\begin{lm}
\[
\frac{(S+P)!}{(S+P-Q)!} = \sum_{I} \left(\begin{array}{c}P\\I\end{array}\right)\left(\begin{array}{c}Q\\I\end{array}\right)\frac{I!S!}{(S-Q+I)!}
\]
\end{lm}
{\bf Proof:} This is a corollary of previous lemma, if one notice that the left hand side of the equation in this lemma is a product with each factor being of the shape of the left hand side of the equation of the previous lemma, and the right hand side is expanded form of a product.\\\\
{\bf Q.E.D.}\\\\
We are actually interested in 
\[
\prod_{j=1}^u\frac{(S+ \sum_{i=1}^j (P_i - Q_i) + Q_j)!}{(S+ \sum_{i=1}^j (P_i - Q_i))!}.
\]
It turns out that we can handle this too.\\
\begin{lm}
\[
\prod_{j=1}^u\frac{(S+ \sum_{i=1}^j (P_i - Q_i) + Q_j)!}{(S+ \sum_{i=1}^j (P_i - Q_i))!}
\]
\[
=\sum_I B_I(L, P, Q) \frac{S!}{(S-(\sum_{j=1}^u Q_j) + I)!}
\]
Here
\[
B_I(L, P, Q) = \sum_{\sum^u_{k=1} I_k = I}\left(\prod^u_{j=1}\left(\begin{array}{c}P_j\\I_j\end{array}\right)\left(\begin{array}{c}\sum^u_{m=j}(Q_m-I_m)+I_j\\I_j\end{array}\right)I_j!\right)
\]
\end{lm}
{\bf Proof:} It is straightforward by suitable induction using previous lemma.\\\\
{\bf Q.E.D.}\\\\
Add in the influence of $\det g$, we get\\
\begin{eqnarray*}
&&\left(\prod_{j=1}^u\frac{(S+ \sum_{i=1}^j (P_i - Q_i) + Q_j)!}{(S+ \sum_{i=1}^j (P_i - Q_i))!}\right)\left(\prod_{j=1}^u \sum_{\tiny{\begin{array}{c}{I_j\subset (1,2, \cdots, r)}\\ {K_j\subset (1,2, \cdots, n)}\end{array}}} \right.\\
&& \left. \det(P^{I_j}_{,K_j})\det(Q^{I_j}_{,K_j}) (L^j)^{e_{I_j}} \frac{(-1)^{|K_j|}}{\prod_{k\in K_j}(s_{,k} + \sum_{i=1}^j(p_{i,k} - q_{i,k}) + q_{j,k})}  \right)\\
&=&\sum_{\tiny{\begin{array}{c}{I_j\subset (1,2, \cdots, r)}\\ {K_j\subset (1,2, \cdots, n)}\\ 1\leq j \leq u \end{array}}}\left(\prod_{j=1}^u \det(P^{I_j}_{,K_j})\det(Q^{I_j}_{,K_j}) (L^j)^{e_{I_j}}(-1)^{|K_j|}\right)\\
&&\left(\prod_{j=1}^u\frac{(S+ \sum_{i=1}^j (P_i - Q_i) + Q_j -e_{K_j})!}{(S+ \sum_{i=1}^j (P_i - Q_i))!}\right)\\
&=&\sum_{\tiny{\begin{array}{c}{I_j\subset (1,2, \cdots, r)}\\ {K_j\subset (1,2, \cdots, n)}\\ 1\leq j \leq u \end{array}}}\left(\prod_{j=1}^u \det(P^{I_j}_{,K_j})\det(Q^{I_j}_{,K_j}) (L^j)^{e_{I_j}}(-1)^{|K_j|}\right)\\
&&\left(\sum_A B_A(\{P_j-K_j\}_{j=1}^u, \{Q_j-K_j\}_{j=1}^u) \frac{S!}{(S-(\sum_{j=1}^u Q_j) + K + A)!}\right)\\
&=& \sum_A \frac{S!}{(S-(\sum_{j=1}^u Q_j) + A)!}\sum_{\tiny{\begin{array}{c}{I_j\subset (1,2, \cdots, r)}\\ {K_j\subset (1,2, \cdots, n)}\\ 1\leq j \leq u \end{array}}}B_{A-K}(\{P_j-K_j\}_{j=1}^u, \{Q_j-K_j\}_{j=1}^u)\\
&&\left(\prod_{j=1}^u \det(P^{I_j}_{,K_j})\det(Q^{I_j}_{,K_j}) (L^j)^{e_{I_j}}(-1)^{|K_j|}\right)\\
&=& \sum_A \hat{B}_A(L, P, Q) \frac{S!}{(S-(\sum_{j=1}^u Q_j) + A)!}
\end{eqnarray*}
Therefore\\
\begin{eqnarray}
\hat{B}_A(L, P, Q)&=& \sum_{\tiny{\begin{array}{c}{I_j\subset (1,2, \cdots, r)}\\ {K_j\subset (1,2, \cdots, n)}\\ 1\leq j \leq u \end{array}}}B_{A-K}(\{P_j-K_j\}_{j=1}^u, \{Q_j-K_j\}_{j=1}^u) \label{ed}\\
&&\left(\prod_{j=1}^u \det(P^{I_j}_{,K_j})\det(Q^{I_j}_{,K_j}) (L^j)^{e_{I_j}}(-1)^{|K_j|}\right)\nonumber
\end{eqnarray}
The relation of $\hat{B}_I(L, P, Q)$ and $D_I(L, P, Q)$ is as following
\begin{equation}
D_I(L, P, Q) = (-1)^l L!\sum_{u=1}^\infty(-1)^u\sum_{\tiny{\sum_{i=1}^u L^i = L }} \frac{1}{\prod^u_{i=1} {(L^i)!}}\hat{B}_I(L, P, Q)
\label{ee}
\end{equation}
Notice that the first few terms of right hand side can be written as
\[
\prod_{j=1}^u\frac{(S+ \sum_{i=1}^j (P_i - Q_i) + Q_j)!}{(S+ \sum_{i=1}^j (P_i - Q_i))!}
\]
\[
=\frac{S!}{(S-\sum_{j=1}^u Q_j)!} +\sum^r_{i=1} \sum^u_{j=1}\left( P_{j,i}\sum^u_{m=j} Q_{m,i} \right) \frac{S!}{(S-\sum_{j=1}^u Q_j+e_i)!} + \cdots
\]
\[
=\frac{S!}{(S- \sum_{i=1}^r l_iQ^i)!} + \sum_{j=1}^r \hat{B}_{1,j}(L, P, Q)\frac{S!}{(S- \sum_{i=1}^r l_iQ^i + e_j)!}+ \cdots 
\]
Here $\hat{B}_{1,j}(L, P, Q)$ is the same as $\hat{B}_{e_j}(L, P, Q)$. From this one can conclude that 
\begin{prop}
\label{eb}
\[
deg_L(D_A(L, P, Q)) = 2a
\]
\end{prop}
{\bf Proof:}\\
Since 
\[
B_A(L, P, Q) = \sum_{\sum^u_{k=1} I_k = A}\left(\prod^u_{j=1}\left(\begin{array}{c}P_j\\I_j\end{array}\right)\left(\begin{array}{c}\sum^u_{m=j}(Q_m-I_m)+I_j\\I_j\end{array}\right)I_j!\right)
\]
Also
\[
deg_L(P_j) = deg_L(\sum^u_{m=j}(Q_m-I_m)+I_j) = 1
\]
and
\[
deg_L\left(\begin{array}{c}P_j\\I_j\end{array}\right) = deg_L\left(\begin{array}{c}\sum^u_{m=j}(Q_m-I_m)+I_j\\I_j\end{array}\right) = |I_j|
\]
We have\\
\[
deg_L(B_A(L, P, Q)) = 2a
\]
Formula (\ref{ed}) gives also
\[
deg_L(\hat{B}_A(L, P, Q)) = 2a
\]
Then by the expression:\\
\[
D_A(L, P, Q) = (-1)^l L!\sum_{u=1}^\infty(-1)^u\sum_{\tiny{\sum_{i=1}^u L^i = L }} \frac{1}{\prod^u_{i=1} {(L^i)!}}\hat{B}_A(L, P, Q)
\]
and lemma \ref{ec}, we get\\
\[
deg_L(D_A(L, P, Q)) = 2a
\]
{\bf Q.E.D.}\\\\
Put above into the expression of $K_m$ as in the last equation of formula (\ref{db}) and use lemma \ref{ea} for variable $S$, we get\\
\[
K_m = \frac{1}{m}\sum_{S-T}\sum_{\tiny{\begin{array}{c} L=(l_1,l_2, \cdots l_r)\not= 0 \\ \sum_{i=1}^r l_i(P^i - Q^i)=T-S\end{array}}}\prod^r_{i=1} (-b_{P^i\bar{Q}^i})^{l_i}\sum_{u=1}^\infty \sum_{\tiny{\sum_{i=1}^u L^i = L }} \frac{1}{\prod^u_{i=1} {(L^i)!}}
\]
\[
(\sqrt{m}z)^{(T-S)}\frac{(-1)^{u-1}}{u}\prod_{j=1}^u\left(\sum_I D_I(L^j, P, Q)(\sqrt{m}z)^{\sum_{k=1}^r l^j_k Q^k-I}(\sqrt{m}\bar{z})^{\sum_{k=1}^r l^j_k Q^k-I}\right)
\]
\[
+ \frac{1}{m}\sum_S (\sqrt{m}z)^S(\sqrt{m}\bar{z})^S \sum_{k=1}^\infty \frac{(-1)^{k-1}}{k}\sum_{\tiny{ \sum_{i=1}^k S_i = S }} \frac{1}{\prod_{i=1}^k S_i!}
\]\\
\[ 
= \frac{1}{m}\sum_{S-T}\sum_{\tiny{\begin{array}{c} L=(l_1,l_2, \cdots l_r)\not= 0 \\ \sum_{i=1}^r l_i(P^i - Q^i)=T-S\end{array}}}\prod^r_{i=1} (-b_{P^i\bar{Q}^i})^{l_i}\sum_{u=1}^\infty \sum_{\tiny{\sum_{i=1}^u L^i = L }} \frac{1}{\prod^u_{i=1} {(L^i)!}}
\]
\[
\frac{(-1)^{u-1}}{u}\prod_{j=1}^u\left(\sum_I D_I(L^j, P, Q)(\sqrt{m}z)^{\sum_{k=1}^r l^j_k P^k-I}(\sqrt{m}\bar{z})^{\sum_{k=1}^r l^j_k Q^k-I}\right)
\]
\[
+ \frac{1}{m}(m|z|^2)
\]\\
\[ 
= \frac{1}{m}\sum_{S-T}\sum_{\tiny{\begin{array}{c} L=(l_1,l_2, \cdots l_r)\not= 0 \\ \sum_{i=1}^r l_i(P^i - Q^i)=T-S\end{array}}}\prod^r_{i=1} (-b_{P^i\bar{Q}^i}(\sqrt{m}z)^{P^i} (\sqrt{m}\bar{z})^{Q^i})^{l_i}\sum_{u=1}^\infty\frac{(-1)^{u-1}}{u}
\]
\[
\sum_{\tiny{\sum_{i=1}^u L^i = L }} \frac{1}{\prod^u_{i=1} {(L^i)!}}\prod_{j=1}^u\left(\sum_I D_I(L^j, P, Q)(\sqrt{m}z)^{-I}(\sqrt{m}\bar{z})^{-I}\right) + |z|^2
\]\\
\[ 
= \frac{1}{m}\sum_{L\not= 0} \prod^r_{i=1} (-b_{P^i\bar{Q}^i}(\sqrt{m}z)^{P^i} (\sqrt{m}\bar{z})^{Q^i})^{l_i}\sum_{u=1}^\infty \frac{(-1)^{u-1}}{u}\sum_{\tiny{\sum_{i=1}^u L^i = L }} \frac{1}{\prod^u_{i=1} {(L^i)!}}
\]
\[
\sum_I \left(\sum_{\tiny{\sum_{j=1}^u I^j = I }} \prod_{j=1}^u D_{I^j}(L^j, P, Q)\right)(\sqrt{m}z)^{-I}(\sqrt{m}\bar{z})^{-I} +|z|^2
\]
\[ 
= |z|^2+\frac{1}{m} \sum_I (\sqrt{m}z)^{-I}(\sqrt{m}\bar{z})^{-I} 
\]
\[
\left( \sum_{L\not= 0} X^L \sum_{u=1}^\infty \frac{(-1)^{u-1}}{u}\sum_{\tiny{\sum_{i=1}^u L^i = L }} \frac{1}{\prod^u_{i=1} {(L^i)!}}E_I(L, P, Q)\right)
\]\\
Where $X=(x^1, x^2, \cdots x^r)$, $x^i = -b_{P^i\bar{Q}^i}(\sqrt{m}z)^{P^i} (\sqrt{m}\bar{z})^{Q^i}$, and
\begin{equation}
E_I(L, P, Q) = \sum_{\tiny{\sum_{j=1}^u I^j = I }} \prod_{j=1}^u D_{I^j}(L^j, P, Q)
\label{ef}
\end{equation}
Recall again that what we want to show is that all the $m^r$ terms with $r>0$ in the expansion of $K_m(z, \bar{z})$ vanish and the terms that do not involve $m$ should equal to $K(z, \bar{z})$. It is not hard to see with the help of lemma \ref{ea} $(ii)$ and from the fact $E_0(L, P, Q) = 1$ that the terms related to $I = 0$ in $K_m$ give exactly the expansion of $K(z, \bar{z})$. Then what we expect clearly is that the $I > 0$ terms only involve negative $m$-power. More precisely we have\\
\begin{eqnarray}
&&K_m = K_m^0 + K'_m \nonumber\\
&&=\left( |z|^2 + \frac{1}{m}\sum_L X^L \sum_{u=1}^\infty \frac{(-1)^{u-1}}{u}\sum_{\tiny{\sum_{i=1}^u L^i = L }} \frac{1}{\prod^u_{i=1} {(L^i)!}}E_0(L, P, Q)\right) \label{eh}
\end{eqnarray} 
\[
+ \frac{1}{m} \sum_{I\not= 0} (\sqrt{m}z)^{-I}(\sqrt{m}\bar{z})^{-I} \left( \sum_L X^L \sum_{u=1}^\infty \frac{(-1)^{u-1}}{u}\sum_{\tiny{\sum_{i=1}^u L^i = L }} \frac{1}{\prod^u_{i=1} {(L^i)!}}E_I(L, P, Q)\right) \nonumber
\]
\begin{prop}
\[
K_m^0 =|z|^2 +\frac{1}{m}\sum_L X^L \sum_{u=1}^\infty \frac{(-1)^{u-1}}{u}\sum_{\tiny{\sum_{i=1}^u L^i = L }} \frac{1}{\prod^u_{i=1} {(L^i)!}}E_0(L, P, Q)= K
\]
\label{ei}
\end{prop}
{\bf Proof:} Clearly\\
\[
E_0(L, P, Q) = D_0(L, P, Q) = (-1)^l L! \sum_{u=1}^\infty(-1)^u\sum_{\tiny{\sum_{i=1}^u L^i = L }} \frac{1}{\prod^u_{i=1} {(L^i)!}}\hat{B}_0(L, P, Q)
\]
By the obvious fact $\hat{B}_0(L, P, Q)=B_0(L, P, Q) = 1$ and formula (\ref{da}), we have\\
\[
E_0(L, P, Q) = 1
\]
This together with lemma \ref{ea}(ii) give us\\
\begin{eqnarray*}
K_m^0 &=&|z|^2 +\frac{1}{m}\sum_L X^L \sum_{u=1}^\infty \frac{(-1)^{u-1}}{u}\sum_{\tiny{\sum_{i=1}^u L^i = L }} \frac{1}{\prod^u_{i=1} {(L^i)!}}\\
&=&|z|^2 + \frac{1}{m}\sum_{i=1}^r x^i\\
&=&|z|^2 + \frac{1}{m}\sum_{i=1}^r (-b_{P^i\bar{Q}^i}(\sqrt{m}z)^{P^i} (\sqrt{m}\bar{z})^{Q^i})\\
&=&|z|^2 + \sum_{i=1}^r c_{P^i\bar{Q}^i}z^{P^i} \bar{z}^{Q^i} = K\\
\end{eqnarray*}
{\bf Q.E.D.}\\\\
It is instructive to analyze the weight of $m$ in various terms. From the definition, one can easily see that $x^i$ all have weight $1$, also notice that $(\sqrt{m}z)^{-I}(\sqrt{m}\bar{z})^{-I}$ has weight $-i$, we should expect
\begin{equation}
F_I(X, P, Q) = \sum_L X^L \sum_{u=1}^\infty \frac{(-1)^{u-1}}{u}\sum_{\tiny{\sum_{i=1}^u L^i = L }} \frac{1}{\prod^u_{i=1} {(L^i)!}}E_I(L, P, Q)
\label{eg}
\end{equation}
being a polynomial of degree less than or equal to $i$. It is easy to deduce from (\ref{eb}) that $\deg_L(E_I(L, P, Q)) = 2i$, which guarantee that $\deg_X(F_I(X, P, Q)) = 2i$ with the help of lemma \ref{ea}. But this falls short to satisfy our need, which requires vanishing of all coefficients of degree larger than $i$. We do not know a coherent way to show the vanishing to all order. But apparently, there is an algorithm to check the vanishing. For any finite order, there only involve finite steps of checking. Although it could be tedious, in principle, it works for checking convergence up to any finite order. Before proceed further, we would like to examine first that what implication we would have if we have already proved $\deg_X(F_J(X, P, Q)) \leq j$ for $j\leq i$.\\
\begin{prop}
Assume that $\deg_X(F_J(X, P, Q)) \leq j$ for $j\leq i$. Then we have\\
\[
\| g_m - g \|_{C^{2i+3}} = O(\frac{1}{m} )
\]
or more precisely\\
\[
\| g_m - g \|_{C^{a,b}} = O(\frac{1}{m} )\ \ \ {\rm for} \ a\leq i+1 \ {\rm or} b\leq i+1
\]
\label{ej}
\end{prop}
{\bf Proof:} Recall from formula (\ref{eh}) and Proposition \ref{ei} that\\
\[
K_m = K + \frac{1}{m} \sum_{I\not= 0} (\sqrt{m}z)^{-I}(\sqrt{m}\bar{z})^{-I} F_I(X, P, Q)
\]
Since in term of $m$, $X$ has weight one. In order for a term involving $X^L$ in $F_J(X, P, Q)$ to have non-negative $m$-power contribution, $l$ has to be bigger than $i$. Its $z$ degree will be $(P\cdot L -j, Q\cdot L -j)$. Our assumptions clearly imply that $F_J(X, P, Q)$ will not contribute non-negative $m$ power if $j\leq i$. Only terms of $F_J(X, P, Q)$ involving $X^L$ for $l>j>i$ will matter. Recall also that for $c_{P\bar{Q}}$ we always have $(p,q)\geq (2,2)$ (with the first term being the curvature). So the corresponding $z$ degree for those terms are\\
\[
(P\cdot L -j, Q\cdot L -j)\geq (2l-j, 2l-j) \geq (i+3, i+3)
\]
These imply the $C^{a,b}$ convergence of \k potential for $a\leq i+2$ or $b\leq i+2$, which in turn imply our proposition.\\\\
{\bf Q.E.D.}\\\\
Especially, for the case $i=0$, namely without any checking, we already have\\
\begin{co}
\[
\| g_m - g \|_{C^3} = O(\frac{1}{m} )
\]
\label{ek}
i.e. induced Bergmann metrics are $C^3$ convergent to the original metric.
\end{co}
\begin{re}
Our previous combinatoric arguments were quite complicated, because we were aiming at the $C^\infty$ convergence. If we only care about $C^3$ convergence, then use our idea, it is not hard to write a simple direct proof of the above result.\\
\end{re}

As an exercise, we will check the case $i=1$.\\\\
Recall that
\[
F_{e_j}(X, P, Q) = \sum_L X^L \sum_{u=1}^\infty \frac{(-1)^{u-1}}{u}\sum_{\tiny{\sum_{i=1}^u L^i = L }} \frac{1}{\prod^u_{i=1} {(L^i)!}}E_{e_j}(L, P, Q)
\]
\[
E_{e_j}(L, P, Q) = \sum^u_{k=1}D_{e_j}(L^k, P, Q)
\]
\[
D_{e_j}(L, P, Q) = (-1)^l L!\sum_{u=1}^\infty(-1)^u\sum_{\tiny{\sum_{i=1}^u L^i = L }} \frac{L!}{\prod^u_{i=1} {(L^i)!}}\hat{B}_{e_j}(L, P, Q)
\]
\[
B_{e_j}(L, P, Q) = \sum^u_{k=1}P_{k,j}\sum^u_{m=k}Q_{m,j}
\]
\[
= \sum^u_{k=1}\sum^r_{i_1=1}l^k_{i_1}P^{i_1}_{,j}\sum^u_{m=k}\sum^r_{i_2=1}l^m_{i_2}Q^{i_2}_{,j}
\]
\[
= \sum^r_{i_1,i_2=1}P^{i_1}_{,j}Q^{i_2}_{,j} \left(\frac{1}{2} \sum^u_{k \not= m}\frac{L^k!}{(L^k-e_{i_1})!}\frac{L^m!}{(L^m-e_{i_2})!} + \sum^u_{k=1}\frac{L^k!}{(L^k-e_{i_1}-e_{i_2})!}\right)
\]
\[
+ \sum^r_{i=1}P^i_{,j}Q^i_{,j}\sum^u_{k=1}\frac{L^k!}{(L^k-e_i)!}
\]
\[
\hat{B}_{e_j}(L, P, Q) = B_{e_j}(L, P, Q) + \sum_{k=1}^{u}\sum_{i=1}^r P^i_{,j}Q^i_{,j}l^k_i(-1)
\]
\[
= \sum^r_{i_1,i_2=1}P^{i_1}_{,j}Q^{i_2}_{,j} \left(\frac{1}{2} \sum^u_{k \not= m}\frac{L^k!}{(L^k-e_{i_1})!}\frac{L^m!}{(L^m-e_{i_2})!} + \sum^u_{k=1}\frac{L^k!}{(L^k-e_{i_1}-e_{i_2})!}\right)
\]
\[
+ \sum^r_{i=1}P^i_{,j}Q^i_{,j}\sum^u_{k=1}\frac{L^k!}{(L^k-e_i)!}-\sum_{k=1}^{u}\sum_{i=1}^r P^i_{,j}Q^i_{,j}l^k_i
\]
\[
= \sum^r_{i_1,i_2=1}P^{i_1}_{,j}Q^{i_2}_{,j} \left(\frac{1}{2} \sum^u_{k \not= m}\frac{L^k!}{(L^k-e_{i_1})!}\frac{L^m!}{(L^m-e_{i_2})!} + \sum^u_{k=1}\frac{L^k!}{(L^k-e_{i_1}-e_{i_2})!}\right)
\]
With the help of lemma \ref{ec}, we get
\[
D_{e_j}(L, P, Q) = \sum^r_{i_1,i_2=1}P^{i_1}_{,j}Q^{i_2}_{,j} \left(\frac{1}{2} (-1)^l \frac{2!L!}{(L-e_{i_1}-e_{i_2})!} + (-1)^{l-1} \frac{L!}{(L-e_{i_1}-e_{i_2})!}\right)
\]
\[
+ \sum^r_{i=1}P^i_{,j}Q^i_{,j} (-1)^l \frac{L!}{(L-e_i)!}- \sum^r_{i=1}P^i_{,j}Q^i_{,j} (-1)^l \frac{L!}{(L-e_i)!}=0
\]
Which is clearly degree one on $L$. So $i=1$ case checking is confirmed. This checking together with Proposition \ref{ej} give us\\
\begin{co}
\[
\| g_m - g \|_{C^5} = O(\frac{1}{m} )
\]
i.e. induced Bergmann metrics are $C^5$ convergent to the original metric.
\end{co}
We also checked the $i=2$ case, which will imply $C^7$ convergence. It is not impossible to check a few steps further. But complicity will increase and one has to stop at some point with this method. The geometric proof presented in the next section will imply all the necessary vanishing in this section.\\\\
\se{Geometric Approach}
In this section, we will give a very simple and complete geometric proof of the $C^\infty$ convergence of Bergmann metrics, which is qualified of being completely different from the combinatoric one presented in the previous section. Therefore, a by-product is that this result imply those unsolved combinatoric problems at the end of the previous section, which still look sort of mysterious.\\\\
Recall that what we want to show is that all the $m^r$ terms with $r>0$ in the expansion of $K_m(z, \bar{z})$ vanish and the terms that do not involve $m$ should equal to $K(z, \bar{z})$. From the expression of $K_m(z, \bar{z})$ in formula (\ref{db}), we can get the following observation, which turns out to be enough to initiate our geometric proof. Consider the expansion of $K_m(z, \bar{z})$ as in formula (\ref{db}) with respect to $z$, let $y_{S\bar{T}}$ be the coefficients of $z^S\bar{z}^T$ term omitting the part that involve negative powers of $m$. Then $y_{S\bar{T}}=f_{S\bar{T}}(m, c_{P\bar{Q}})$ is a function of $m$ and some $c_{P\bar{Q}}$.\\\\
\begin{prop}
\label{fa}
$y_{S\bar{T}}=f_{S\bar{T}}(m, c_{P\bar{Q}})$ is a polynomial of m and $c_{P\bar{Q}}$ for $p+q \leq s+t$ with degree of each variable bounded from above with respect to $S$, $T$.
\end{prop}
{\bf Proof:}
It is clear, if one observe that $z^S\bar{z}^T$ is accompanied by the factor $m^\frac{s+t}{2}$, and all $b_{P\bar{Q}}$ contain a negative $m$-power of $m^{1-\frac{p+q}{2}}$. So $\deg_m f_{S\bar{T}} \leq \frac{s+t}{2}$ and $\deg_{b_{P\bar{Q}}}$ large terms will have negative $m$-power and will be omitted.\\
{\bf Q.E.D.}\\\\
Let $\{ f_m \}_{m=1}^\infty$ be a sequence of real functions, such that $\lim_{m\rightarrow \infty} f_m = f_0$. The idea of our geometric prove can be best illustrated in the following lemma.\\
\begin{lm}
\label{fb}
For $r\leq 0$, assume that $\lim_{m\rightarrow \infty} f_m = f_0$ ,and $\lim_{m\rightarrow \infty} m^r f'_m = g$,\\
(i) if $r<0$, then $g = 0$,\\
(ii) if $r=0$, then $g = f'_0$
\end{lm}
{\bf Proof:}
Take a test function $\varphi$, then $\lim_{m\rightarrow \infty} m^r f'_m = g$ implies:
\[
\int g \varphi = \lim_{m\rightarrow \infty} \int m^r f'_m \varphi = \lim_{m\rightarrow \infty} - m^r\int f_m \varphi' 
\]
\[
= \lim_{m\rightarrow \infty} (- m^r)\lim_{m\rightarrow \infty}\int f_m \varphi' = (\lim_{m\rightarrow \infty} - m^r) \int f_0 \varphi'
\]
\[
= \left\{
\begin{array}{cc}
0, \ & \ r < 0\\
\int f'_0 \varphi, \ & \ r = 0
\end{array}
\right.
\]
Therefore
\[
g= \left\{
\begin{array}{cc}
0 \ & \ r < 0\\
f'_0 \ & \ r = 0
\end{array}
\right.
\]
{\bf Q.E.D.}\\\\
For a point $m_0 \in M$, we would like to discuss the convergence on the canonical coordinate chart: $w = \phi_{m_0} : U_{m_0} \longrightarrow {\bf C^n}$. Recall that for $\phi_{m_0}^{-1}(w) \in U_{m_0}$, we can choose a smooth family of canonical coordinates $z = \phi_w : U_w \longrightarrow {\bf C^n}$ with the convention that $\phi_0 = \phi_{m_0}$, which amounts to pick a smooth section of the principle bundle of orthonormal frames over $U_{m_0}$ passing through the canonical frame at $m_0$.\\\\
To prove the convergence of $\{ g_m \}$ to $g$, we will use the induction on the order of convergence. Assume that $\{ g_m \}$ converge to $g$ up to $C^{l-1}$, we need to show that for any $L$ with respect to $l$ ($L=(P,Q)$, $p+q=l$), $\partial_L g_m$ converge to $\partial_L g$. Under coordinate $w$, let $g(w)= g_{ij}dw^id\bar{w}^j$, and $g_m(w)= g_{m,ij}dw^id\bar{w}^j$, written in local coordinate, we need to show that $\partial_L g_{m,ij}$ converge to $\partial_L g_{ij}$. Notice that 
\[
\partial_L g_{ij}(w) = \partial_{L+(e_i,e_j)} K^w (w) = \sum_M \frac{{\partial z}^M}{{\partial w}^{L+(e_i,e_j)}}\partial_M K^w (z) + h_L
\]
\[
\partial_L g_{m,ij}(w) = \partial_{L+(e_i,e_j)} K^w_m (w) = \sum_M \frac{{\partial z}^M}{{\partial w}^{L+(e_i,e_j)}}\partial_M K^w_m (z) + h_{m,L}
\]
Here $\deg M = (p+1,q+1)$, and $\frac{{\partial z}^M}{{\partial w}^{L+(e_i,e_j)}}$ indicate a sum of fractions with nominators being permutation of factors, which are exactly the factors one get for the highest degree derivatives after using chain rule. Notice that at $w=0$, $\frac{{\partial z}^M}{{\partial w}^{L+(e_i,e_j)}}= \delta ^M_{L+(e_i,e_j)}$. $h_L$ and $h_{m,L}$ indicate the parts that contain derivatives of lower order. by our inductive assumption, $\{h_{m,L}\}$ converge to $h_L$. We only need to deal with the highest order terms. According to Proposition \ref{fa}, the part of $\partial_L g_{m,ij}(w)$ that do not involve negative powers of $m$ (denote by $g_L(w, m, c_{P\bar{Q}}(w))$) is a polynomial on $m$ and some $c_{P\bar{Q}}(w)$ of bounded degree. Let $\deg_m(g_L)=d_L$, then
\[
g_L(w, m, c_{P\bar{Q}}(w)) = \sum_{i=1}^{d_L} g_{L,i}(w) m^i
\]
we can conclude that
\[
\lim_{m\rightarrow \infty} m^{-d_L} \partial_L g_{m,ij}(w) = \lim_{m\rightarrow \infty} m^{-d_L}g_L(w, m, c_{P\bar{Q}}(w)) = g_{L,d_L}(w)
\]
According to lemma \ref{fb}, we have 
\[
\lim_{m\rightarrow \infty} m^{-d_L} \partial_L g_{m,ij}(w) = g_{L,d_L}(w) =0
\]
unless $d_L = 0$. Therefore $\deg_m(g_L)= 0$ and use induction assumption and lemma \ref{fb} again, we get
\[
\lim_{m\rightarrow \infty} \partial_L g_{m,ij}(w) = \partial_L g_{ij}(w)
\]
Now by the induction, we have therefore proved the $C^\infty$ convergence of the Bergmann metrics, as stated in the following theorem, modulo the starting point of the induction, the $C^0$ convergence of the Bergmann metrics, which we can refer to corollary \ref{ek} in previous chapter or Tian's proof, although it is not hard to give it a direct proof.\\
\begin{th}
For any $l\in {\bf Z_+}$ 
\[
\| g_m - g \|_{C^l} = O_l(\frac{1}{m} )
\]
i.e. induced Bergmann metrics are $C^\infty$ convergent to the original metric.
\end{th}
The subindex $l$ in the right hand side indicate its dependence on $l$.\\\\
A straight forward corollary of this theorem is those combinatoric vanishing problems in the end of the previous section, which we was not able to solve directly.\\
\begin{co}
\[
\deg_X(F_I(X, P, Q)) = i
\]
Here $F_I(X, P, Q)$ is defined in the previous section. 
\end{co}
\begin{re}
Anybody, who have ventured through last section, must have realized how simple the geometric proof is. Here again, we witness the power of analysis. Nevertheless this proof is not as straight forward as the combinatoric one. Does not really help to uncover the combinatoric mystery in the previous section. For this reason, we still like to present the combinatoric approach toward the problem, which seems quite interesting in its own.
\end{re}
\se{Application}
In this section we will discuss an interesting application of the convergence of Bergmann metric. For $(M, \omega)$ a \k manifold with the associated positive Hermitian line bundle $(L, h)$, where $c_1(L, h)=\omega$, in arithmetic geometry, one usually need some quantitive result related to the Hermitian metric. A typical result is the following:\\
\begin{th}
For any $\epsilon > 0$, there exist a $n_0$, such that for any $p \in M$ and $n> n_0$, there is a section $s \in \Gamma(L^{\otimes n})$ such that
\[
\|s\|_{sup,h}\leq e^{n\epsilon}\|s\|_h(p)
\]\\
\end{th}
It is then interesting to see if the theorem will still be true with $\epsilon = 0$. We will show as an application of our techniques and results that this sharp version and some generalization of it is still true.\\
\begin{th}
There exist a $n_0$ depends only on geometry of $(L, h)$, such that for any $p \in M$ and $n> n_0$, there is a section $s \in \Gamma(L^{\otimes n})$ such that
\[
\|s\|_{sup,h}= \|s\|_h(p)
\]\\
\end{th}
{\bf Proof:} For simplicity we will first deal with the case that $L$ is very ample. For a point $p\in M$, we will construct a section of $L^m$, for which $\|s\|_h$ achieves its maxima at $p$. (Actually we will show that $p$ is the only maximal point.) Let $z = (z^1, z^2, \cdots, z^n)$ be the local K-coordinate of $(M, g)$ with respect to $p$. we will prove the following three statements which together will imply the theorem.\\\\
(i) For certain big $C_2>0$ and $|z|\geq \frac{C_2}{m}$, when $m$ is large enough we have
\[
\|s\|_h(0) > \|s\|_h(z)
\]
(ii) For certain small $C_1>0$ and $|z|\leq \frac{C_1}{m}$, when $m$ is large enough we have
\[
\|s\|_h(0) > \|s\|_h(z)
\]
(iii) For certain small $C_2>C_1>0$ and $\frac{C_1}{m}\leq |z|\leq \frac{C_2}{m}$, when $m$ is large enough we have
\[
\|s\|_h(0) > \|s\|_h(z)
\]\\
Since we assume that $L$ is very ample, we have an embedding:\\
\[
i: M \longrightarrow {\bf CP}^l
\]
Where ${\bf CP}^l$ has coordinate $W = [W_0, W_1, \cdots, W_l]$, $i^* O(1) = L$. $W_i$ can be viewed as sections of $O(1)$, we will also use them to denote pull back sections of $L$. We assume that $i(p) = [1, 0,  \cdots, 0]$.\\\\
An important observation is that for $O(1)$ on ${\bf CP}^l$, it is very easy to find ``peak'' section needed. $W_0$ will work. For $O(m)$, we may take $(W_0)^m$. One would hope that $(W_0)^m$ will also work for $L^m$. But $g$ is not the pull back of the Fubini-Study metric. We need to modify it a little. The key observation here is that the term that spoil the ``peak'' property is the first order term. So we will need to adjust by some first order term.\\\\
Choose $W$ suitably such that for $w= (w^1, \cdots, w^n)$, where $w^i= \frac{W_i}{W_0}$, we have\\
\[
\left. \frac{\partial w}{\partial z}\right|_{z=0} = id
\]
Also let $h_0$ denotes the Hermitian metric on $L=i^* O(1)$ induced from usual Hermitian metric on $O(1)$. For peak section $e_p$ of $(L, h)$ at point $p\in M$, we have\\
\[
\left\|e_p^m\right\|_h = e^{-mK}, \ {\rm and}\ \left\|e_p^m\right\|_{h_0} = 1 + \sum_{S,T}\frac{G^{T,\bar{S}}}{G^{0,\bar{0}}} z^S \bar{z}^T
\]
Then for any section $s$ of $L^m$, we have\\
\[
\left\|s\right\|_h = f(z)\left\|s\right\|_{h_0}
\]
where 
\[
f(z) = \frac{e^{-mK}}{1 + \sum_{S,T}\frac{G^{T,\bar{S}}}{G^{0,\bar{0}}} z^S \bar{z}^T} = e^{-m(K-K_m)}
\]
By our theorem $|K-K_m| = O(\frac{1}{m})$, so we have\\
\[
f(z) = 1 - \sum_i\left(\frac{G^{0,\bar{e_i}}}{G^{0,\bar{0}}} z^i + \frac{G^{e_i,\bar{0}}}{G^{0,\bar{0}}}\bar{z}^i\right) + \sum_{S,T} f_{S,\bar{T}} z^S \bar{z}^T
\]
where $f_{S,\bar{T}} = O(1)$ and $s+t\geq 2$. AN important fact is\\
\begin{equation}
G^{0,\bar{e_i}}=O(\frac{1}{m})
\label{ga}
\end{equation}
Although apriorily we only have $G^{0,\bar{e_i}}=O(1)$. Recall the formula\\
\[
G^{0,\bar{e_i}} = m^{\frac{1}{2}} \sum_{\tiny{\begin{array}{c} L=(l_1,l_2, \cdots l_r) \\ \sum_{i=1}^r l_i(P^i - Q^i)=e_i\end{array}}}\frac{\prod^r_{i=1} (-b_{P^i\bar{Q}^i})^{l_i}}{L!}F(L,P,Q)
\]
also that\\
\[
b_{P\bar{Q}} = -c_{P\bar{Q}} m^{1-\frac{p+q}{2}}
\]
Notice that $c_{P\bar{Q}}$ do not depend on $m$ and $p,q\geq 2$, so all the terms of $G^{0,\bar{e_i}}$ will be of order at least $O(\frac{1}{m})$ unless $L=e_j$ and $p^j=q^j =2$ which is impossible to satisfy $l_i(P^i - Q^i)=e_i$. So we have that estimate (\ref{ga}) is correct. Let\\
\[
s_0 = (W_0)^m,\ s_1= \left(\sum_i \frac{G^{0,\bar{e_i}}}{G^{0,\bar{0}}} W_i\right) (W_0)^{m-1}\ {\rm and}\ s=s_0+s_1
\]
Since\\
\[
|K-K_m|=O(\frac{1}{m}) \ {\rm and}\ e^{-m(K-K_m)} = O(1)
\]
We always have that
\[
\|s\|_h \leq C \|s\|_{h_0}
\]
In case (i) we have\\
\[
\|s\|_h(0) =\|s_0\|_h(0) = \|s_0\|_{h_0}(0) =1
\]
also\\
\[
\|s\|_h(z) \leq \|s_0\|_h(z) + \|s_1\|_h(z)
\]
\[
\leq C(\|s_0\|_{h_0}(z) +\|s_1\|_{h_0}(z))
\]
By estimate (\ref{ga}) we have
\[
\|s_1\|_{h_0}(z) = O(\frac{1}{m}) 
\]
It is easy to see when $C_2$ is large enough, $\|s_0\|_{h_0}(z)$ can be as small as we want. Together we have\\
\[
\|s\|_h(0) \geq \|s\|_h(z)
\]\\
For case (iii), notice
\[
e^{-mK} -e^{-mK_m} = e^{-mK_m}(e^{-m(K-K_m)} -1)
\]
\[
m|K-K_m| = -\sum_i\left(\frac{G^{0,\bar{e_i}}}{G^{0,\bar{0}}} z^i + \frac{G^{e_i,\bar{0}}}{G^{0,\bar{0}}}\bar{z}^i\right) + O(|z|^2)
\]
Since we have $|z|^2<\frac{C_2}{m}$, we have\\
\[
m|K-K_m| = O(\frac{1}{m})
\]
So in this case we have\\
\[
\|s\|_h -\|s\|_{h_0} = O(\frac{1}{m})\|s\|_{h_0}
\]
For any $z$ such that $\frac{C_1}{m}\leq |z|^2\leq \frac{C_2}{m}$, we have\\
\[
\|s\|_h(z)= \|s_0\|_{h_0}(z)+(\|s\|_h(z)-\|s\|_{h_0}(z))+(\|s\|_{h_0}(z) - \|s_0\|_{h_0}(z))
\]
\[
\leq \|s_0\|_{h_0}(z) + O(\frac{1}{m})\|s\|_{h_0} (z) + \|s_1\|_{h_0}(z)
\]
\[
\leq (1+O(\frac{1}{m}))\|s_0\|_{h_0}(z) + \|s_1\|_{h_0}(z)
\]
It is easy to show that for $z$ in our range\\
\[
\|s_0\|_{h_0}(z) \leq C_3 <1
\]
where $C_3$ depends on $C_1$. When $m$ is large enough, we have\\
\[
\|s\|_h(z) \leq (1+O(\frac{1}{m}))C_3 + \frac{1}{m}|z| < 1 = \|s\|_h(0)
\]
For case (ii), it is easy to show that\\
\[
\|s\|_h = 1-Cm|z|^2 +O(m^2|z|^4 + m|z|^3)
\]
which clearly achieves its maximum at $z=0$ for $|z|^2\leq \frac{C_1}{m}$ with $C_1$ small enough.\\\\
Now all the three cases are proved.\\\\
{\bf Q.E.D.}\\\\
{\bf Acknowledgment}: I would like to thank Prof. M. Kontsevich for very stimulating discussion on mirror symmetry and Prof. Shouwu Zhang for proposing the problem in the last section. I also thank Prof. S. T. Yau, Prof. M. Weber and Prof. K. Koehler for pointing out the references leading to [B]. The work is done while my stay in IHES of France and Ruhr University at Bochum of Germany. I would like to thank both institutes for Hospitality.\\\\\\
\noindent 
{\bf Reference:\\}\\

[B] S. Bochner, {\em Curvature in Hermitian metric}, Bull. Amer. Math. Soc., {\bf 53} (1947), 179-195\\

[T] G. Tian, {\em On a set of polarized \k metrics on algebraic manifolds}, J. Differential Geometry  {\bf 32} (1990), 99-130.\\

\end{document}